\documentclass[12pt,a4paper]{article}

\usepackage[dvips]{graphics}
\usepackage[dvips]{epsfig}

\newtheorem{lemma}{Lemma}
\newtheorem{thm}[lemma]{Theorem}
\newtheorem{cor}[lemma]{Corollary}
\newtheorem{prop}[lemma]{Proposition}

\newtheorem{rem}{Remark}

\newcommand{\dimo}{\noindent \emph{Proof. }}
\newcommand{\qed}{\\ \vspace{-0.3truecm} \rightline{$\Box$ \ \ \ \ \ \ \ \ \ \ \ \ \ \ \ }\\}

\newcommand{\e}{\varepsilon}

\usepackage{latexsym}
\usepackage{amsfonts}
\usepackage{amssymb}
\usepackage{amsmath}
\usepackage[english]{babel}

\begin{document}
\title{PL $4$-MANIFOLDS \\ ADMITTING SIMPLE CRYSTALLIZATIONS: \\ FRAMED LINKS AND REGULAR GENUS}

\author{Maria Rita CASALI  - Paola CRISTOFORI - Carlo GAGLIARDI \\
Dipartimento di Scienze Fisiche, Informatiche e  Matematiche \\
       Universit\`a di Modena e Reggio Emilia \\
       Via Campi 213 B \\  I-41125 MODENA (Italy)}

\maketitle

\abstract {{\it Simple crystallizations} are edge-coloured graphs representing PL 4-manifolds with the property that the 1-skeleton of the associated triangulation equals the 1-skeleton of a 4-simplex.
In the present paper, we prove that any (simply-connected) PL $4$-manifold $M$ admitting a simple crystallization admits a {\it special handlebody decomposition}, too; equivalently, $M$ may be represented by a framed link yielding $\mathbb S^3$, with exactly $\beta_2(M)$ components ($\beta_2(M)$ being the second Betti number of $M$).
As a consequence, the {\it regular genus} of $M$ is proved to be the double of $\beta_2(M)$.

Moreover, the characterization of any such PL $4$-manifold by $k(M)= 3 \beta_2(M)$, where $k(M)$ is the {\it gem-complexity} of $M$ (i.e. the non-negative number $p-1$, $2p$ being the minimum order of a crystallization of $M$) implies that both PL invariants gem-complexity and regular genus turn out to be additive within the class of all PL $4$-manifolds admitting simple crystallizations (in particular: within the class of all  ``standard'' simply-connected PL 4-manifolds).
} \endabstract



\bigskip

\section{\hskip -0.7cm . Introduction and main results}

For any PL $n$-manifold $M^n$, it is known the existence of a {\it contracted triangulation}, i.e. a {\it pseudocomplex}\footnote{Roughly speaking, pseudocomplexes extend the notion of simplicial complexes, since a set of vertices may determine more than one face. In the literature, a similar notion is also given by the term {\it simplicial poset}.} triangulating $M^n$, whose $0$-skeleton consists of exactly $n+1$ vertices.

\medskip

Note that contracted triangulations of a PL $n$-manifold $M$  may be seen as an intermediate notion, between {\it simplicial complexes} (where a pair of distinct simplices can intersect in at most one face) and  {\it singular triangulations}  (where it is required that only the interior of the cells are open simplices, and usually a single 0-simplex is present).

What makes contracted triangulations particularly user-friendly is the possibility of representing them by means of their dual graphs, which turn out to be a special kind of edge-coloured graphs, called {\it crystallizations}.
Hence, results in PL-topology can be obtained within {\it crystallization theory}\footnote{See \cite{[FGG]}, \cite{survey} and \cite{Casali-Oberwolfach}  for surveys on crystallizations by the founding group of the theory, and \cite{[BaD]}, \cite{[BaS]}, \cite{[Sw]} for some recent results by different authors which are contributing to its development.} via  combinatorial tools.

As an example, catalogues of PL 3-and 4-manifolds have been automatically generated and classified for increasing order of the associated crystallizations (see \cite{[CC$_2$]} and \cite{[CC$_2014$]}, together with their references) and the classifying algorithm for PL manifolds in dimension $4$ appears to be a promising approach to the problem of detecting different PL-structures on the same topological 4-manifold (\cite[Section 5]{[CC$_2014$]}).
Furthermore, the representation by crystallizations has allowed the definition of graph-defined PL-invariants: one of them, the {\it regular genus} (\cite{[G]}), has yielded classification theorems of particular
significance in dimension 4 and 5 (see, for example,  \cite{[CG]}, \cite{Casali} and \cite{Casali-Malagoli}).

In this context, the needs coming from the analysis of the catalogues on one hand and from the computation of the regular genus on the other hand make of considerable interest the identification of crystallizations of a PL 4-manifold that are minimal with respect to the number of vertices or with respect to the regular genus.
In this paper we show that {\it simple crystallizations}, recently introduced in \cite{[BaS]}, meet both requirements, and we study the properties of the PL $4$-manifolds admitting such crystallizations.
A crystallization is called {\it simple} if the 1-skeleton of the associated contracted triangulation equals the 1-skeleton of a 4-simplex.

\bigskip

The key step of the analysis is the proof that any contracted triangulation associated to a simple crystallization induces a particular type of handle decomposition lacking in 1-handles and 3-handles (i.e. a {\it special handlebody decomposition}, according to \cite[p. 59]{[M]}): see Proposition \ref{handle-decomposition}.\footnote{Note that the existence of a special handlebody decomposition is related to Kirby problem n. 50, and is of particular interest with regard to exotic PL 4-manifolds: see, for example, \cite{[Ak1]} and  \cite{[Ak2]}.}
Equivalently, if $M$ is the represented PL 4-manifold, a framed link $L$  with $\beta_2(M)$ components exists, so that the 3-dimensional Dehn's surgery on $L$ yields $\mathbb S^3$, while in dimension $4$ the same framed link $L$ (without dotted circles) identifies both the bounded PL 4-manifold $M-{\mathbb D^4}$ and the closed PL 4-manifold $M$ itself.
Hence:

\begin{thm} \label{framed-link}
\ \
If $M$ admits simple crystallizations, then $M$ is represented by a (not dotted) framed link  with $\beta_2(M)$ components.
\end{thm}

As a consequence, various combinatorial properties of simple crystallizations are obtained  (Proposition \ref{proprieta_cryst_semplici}  and Proposition \ref{minimality}), which
allow to prove that the regular genus of PL 4-manifolds admitting simple crystallizations  is twice their Betti number.
Further, such PL 4-manifolds turn out to be characterized by a nice relation between their {\it gem-complexity}  (i.e. the PL invariant $k(M)=p-1$, where $2p$ is the minimum order of a crystallization of $M$: see Definition 1) and their second Betti number.

\begin{thm}  \label{main}
\ \
If $M$ admits simple crystallizations, then \  $ \mathcal G(M) \, = \, 2 \beta_2(M).$
\par \noindent Moreover, a closed simply-connected PL 4-manifold $M$ admits simple crystallizations  if and only if  \ $k(M) \, = \, 3 \beta_2(M).$
\end{thm}

In virtue of Theorem \ref{main},  both the invariants gem-complexity and regular genus turn out to be additive with respect to connected sum within the class of all PL $4$-manifolds admitting simple crystallizations (in particular: within the class of all  ``standard'' simply-connected PL 4-manifolds): see Proposition \ref{additivity}.

Note that subadditivity of both gem-complexity and regular genus is known to hold in any dimension, while the additivity of regular genus has been conjectured, but the problem is still open (with the exception of the low-dimensional cases).
In particular, in dimension four, additivity of regular genus - at least in the simply-connected case - would imply the 4-dimensional Smooth Poincar\'e Conjecture: see \cite[Remark 1]{[FG$_2$]} (resp. \cite[Remark 2]{[CC$_2014$]}).
From this viewpoint, our result about additivity for 4-manifolds admitting simple crystallizations appears to be significant, in connection with the problem of the existence of simple crystallizations for a given simply-connected PL $4$-manifold (see Proposition \ref{exotic-simple} for some particular families  yielding a negative answer and for relationships with 4-dimensional crystallization catalogues).
Finally, we point out that simple crystallizations may be useful in order to prove algorithmically  the PL-equivalence of different triangulations of the same (simply-connected) topological 4-manifold: for example, in \cite[Section 1]{[BaS]} and \cite[Section 1]{[CC$_2014$]}, ongoing attempts to prove via simple crystallizations the conjecture of \cite{[SK]}
concerning the $K^3$-surface are described.

\bigskip

\section{\hskip -0.7cm . Basic notions on coloured triangulations of PL manifolds}

Edge-coloured graphs are a representation tool for the whole class of piecewise linear (PL) manifolds, without restrictions about dimension, connectedness, orientability or boundary properties. In the present work, however, we will deal only with closed, connected and orientable PL-manifolds of dimension $n=4$; hence, we will briefly review basic notions and results of the theory with respect to this particular case.

A {\it 5-coloured graph (without boundary)} is a pair
$(\Gamma,\gamma)$, where $\Gamma= (V(\Gamma),$ $E(\Gamma))$ is a
regular multigraph (i.e. it may include multiple edges, but no
loop) of degree five and $\gamma : E(\Gamma) \to
\Delta_4=\{0,1,2,3,4\}$ is a proper edge-coloration (i.e. it is
injective when restricted to the set of edges incident to any vertex of $\Gamma$).

\smallskip

The elements of the set $\Delta_4$ are called the {\it colours} of
$\Gamma$; thus, for every $i\in \Delta_4$, an {\it $i$-coloured
edge} is an element $e \in E(\Gamma)$ such that $\gamma(e)=i.$ For
every $i,j,k \in \Delta_4$ let $\Gamma_{\hat \imath}$ (resp. $\Gamma_{ijk}$) (resp.
$\Gamma_{ij}$) be the subgraph obtained from $(\Gamma, \gamma)$ by deleting all the
edges of colour $i$ (resp. $c\in \Delta_4-\{i,j,k\}$) (resp. $c\in \Delta_4-\{i,j\}$).
The connected components of $\Gamma_{\hat\imath}$ (resp. $\Gamma_{ijk}$) (resp. $\Gamma_{ij}$)
are called {\it ${\hat\imath}$-residues}  (resp. {\it $\{i,j,k\}$-coloured residues}) (resp. {\it $\{i,j\}$-coloured cycles}) of $\Gamma,$ and their number is denoted by  $g_{\hat\imath}$ (resp. $g_{ijk}$) (resp. $g_{ij}$).
A 5-coloured graph $(\Gamma, \gamma)$ is called {\it contracted} iff, for each $i\in \Delta_4$, the subgraph
$\Gamma_{\hat\imath}$ is connected (i.e. iff $g_{\hat\imath}=1$ $\forall i\in \Delta_4$).

Every 5-coloured graph $(\Gamma, \gamma)$ may be thought of as the
combinatorial visualization of a $4$-dimensional labelled
pseudocomplex $K(\Gamma)$, which is constructed
according to the following instructions:
\begin{itemize}
\item{}
\  for each vertex $v\in V(\Gamma)$, take a $4$-simplex $\sigma(v)$, with vertices labelled $0,1,2,3,4$;
\item{}
\ for each $j$-coloured edge between $v$ and $w$ ($v,w\in V(\Gamma)$), identify the 3-dimensional  faces of $\sigma(v)$ and $\sigma(w)$ opposite to the vertex labelled $j$, so that equally labelled vertices coincide.
\end{itemize}

In case $K(\Gamma)$ triangulates a (closed) PL $4$-manifold $M$, then it is called a  {\it coloured triangulation} of $M$ and  $(\Gamma,\gamma)$  a {\it gem} (gem = \underbar graph  \underbar encoded \underbar manifold) {\it representing} $M$.

\smallskip

The construction of $K(\Gamma)$ directly ensures
that, if $(\Gamma,\gamma)$ is an order $2p$  gem of $M$, then:
\begin{itemize}
\item[(a)] $M$ is orientable iff $\Gamma$ is bipartite;
\item[(b)] there is a bijection between $i$-labelled vertices (resp. $1$-simplices whose vertices are labelled $\Delta_4 - \{i,j,k\}$) (resp. $2$-simplices whose vertices are labelled $\Delta_4 -\{i,j\}$) of
$K(\Gamma)$  and ${\hat\imath}$-residues  (resp. $\{i,j,k\}$-coloured residues) (resp. $\{i,j\}$-coloured cycles) of $\Gamma;$
\item[(c)] $\chi (|K(\Gamma)|) = - 3 p + \sum_{i,j} g_{ij} - \sum_{i,j,k} g_{ijk} + \sum_i g_{\hat \imath};$
\item[(d)]  $2 g_{ijk}  =  g_{ij} + g_{ik} + g_{jk} - p$ for each triple $(i,j,k)\in \Delta_4$.
  \end{itemize}

Finally, a gem representing a (closed) PL $4$-manifold $M$ is  a {\it crystallization} of $M$ if it is also a contracted graph; by the above property (b), this is equivalent to requiring that the associated pseudocomplex
$K(\Gamma)$ contains exactly five vertices (one for each label $i\in \Delta_4$).
Pezzana's Theorem and its subsequent improvements (\cite{[FGG]}) prove that every PL-manifold admits a crystallization.

\bigskip

As already recalled, catalogues of PL manifolds have been obtained both in dimension three (see \cite{[L]}, \cite{[CC$_1$]} and \cite{[CC$_2$]} for the 3-dimensional orientable case and \cite{[C$_2$]}, \cite{[C$_4$]} and \cite{[BCrG$_1$]}  for the non-orientable one) and four (\cite{[CC$_2014$]}).  They are constructed with respect to a suitable graph-defined PL invariant, which measures how ``complicated'' is the representing combinatorial object:\footnote{The approach is similar to Matveev's for 3-dimensional censuses, where the 3-manifolds are listed with respect to  the minimum number of {\it true vertices} in their special spines. See \cite{[Mat]} for details about {\it Matveev complexity},  and  \cite{[CC$_1$]} and \cite{[CC$_3$]} for relationships with gem-complexity.}
\smallskip

\par \noindent
\textbf{Definition 1.} \ Given a PL $n$-manifold $M^n$, its \textit{gem-complexity} is the non-negative integer $k(M^n)= p - 1$, where $2p$ is the minimum order of a crystallization of $M^n$.

\medskip

Note that, as proved for arbitrary dimension $n \ge 3$ in \cite[Proposition 7]{[CC$_2014$]}, if  $M$ is assumed to be a {\it handle-free} PL 4-manifold (i.e.: if it admits neither the orientable nor the non-orientable $\mathbb S^3$-bundle over $\mathbb S^1$ as a connected summand), then $k(M)=p-1,$ where $2p$ is the order of a crystallization of $M$ lacking in {\it 2-dipoles} (i.e. pairs of parallel edges coloured by $\Delta_4-\{i,j,k\}$, whose end-points belong to different $\{i,j,k\}$-residues) and {\it $\rho$-pairs} (i.e.: pairs of distinct $i$-coloured edges both belonging to at least three common bicoloured cycles).

Crystallizations with these properties are called {\it rigid dipole-free crystallizations}; they are exactly the elements considered in the existing crystallization catalogues in dimension four.\footnote{A slightly modified definition of {\it rigidity} is required in 3-dimensional crystallization catalogues.}

\bigskip

As mentioned in Section 1, some of the most interesting results of crystallization theory are related to a graph-based invariant for PL $n$-manifolds, called {\it regular genus} and introduced in \cite{[G]}\footnote{See, for example, \cite{[C$_1992$]}, \cite{[Casali_UMI]} and \cite{Casali-Malagoli} for 4-dimensional results, \cite{[CG]}  and \cite{Casali} for 5-dimensional ones.}.
It extends to arbitrary dimension the classical notion of Heegaard genus of a 3-manifold
and relies on the existence of a particular type of embedding into a surface for graphs representing manifolds of arbitrary dimension.

As far as the 4-dimensional case is concerned, it is well-known that, if $(\Gamma, \gamma)$ is an order $2p$ crystallization of an orientable\footnote{Analogous results and definitions exist in the non-orientable case, too, as well as in general dimension  (see \cite{[G]}); for the purpose of the present work, however, the attention may be restricted to bipartite graphs representing orientable PL 4-manifolds.} PL 4-manifold $M,$ then for every  cyclic permutation $\e= (\e_o, \e_1, \e_2, \e_3, \e_4=4)$ of $\Delta_4$  there exists a so-called {\it regular embedding}\footnote{By short, it is a cellular embedding whose regions are bounded by the images of $\{\e_i,\e_{i+1}\}$-coloured cycles, for each $i \in \mathbb Z_5$.}
$i_{\e}\ :\ |\Gamma| \to F_{\e}$, where $F_{\e}$ is a closed orientable surface whose  genus - denoted by $\rho_\e(\Gamma)$ - may be directly computed by the following formula (see \cite{[G]} for details):
\begin{equation}\label{genere}
\sum_{i\in \mathbb Z_{5}} g_{\e_{i}\e_{i+1}} - 3 p  = 2 - 2 \rho_\e(\Gamma).
\end{equation}

\par \noindent
\textbf{Definition 2.} \ The {\it regular genus} of a bipartite 5-coloured graph $\Gamma$  is defined as the minimum genus of a surface into which $\Gamma$ regularly embeds:
$$\rho(\Gamma) = \min_\e \{\rho_\e(\Gamma)\};$$
the {\it regular genus} of a PL 4-manifold $M$ is defined as the minimum regular genus of a crystallization of $M$:
$$\mathcal G(M) = \min \{\rho(\Gamma) \ / \ (\Gamma, \gamma) \ \text{crystallization of} \  M\}.$$

\bigskip

For the purpose of the present paper, it is worthwhile to note that, if the PL 4-manifold $M$ is assumed to be simply-connected, the following relation involving the regular genus and the second Betti number $\beta_2(M)$ of $M$ always holds (see equality (5) in \cite{[CC$_2014$]}, or \cite[Proposition 2]{[Cav$_1992$]}):
\begin{equation}\label{genere-beta_2}
 \beta_2(M) \le \left[ \frac {\mathcal G(M)}2 \right],
\end{equation}
where $[x]$ denotes the integer part of $x$.

\bigskip

\section{\hskip -0.7cm . 4-manifolds admitting simple crystallizations}

In \cite{[BaS]} the notion of {\it simple crystallization} of a (simply-connected) PL $4$-manifold is introduced:

\smallskip

\par \noindent
\textbf{Definition 3.} \   A $4$-dimensional pseudocomplex $K$ triangulating a PL 4-manifold $M$ is said to be {\it simple} if any pair of vertices belongs to at most one $1$-simplex.   A {\it simple crystallization} of a PL 4-manifold $M$ is a crystallization $(\Gamma, \gamma)$ of $M$, whose associated pseudocomplex $K(\Gamma)$ is simple.
\medskip

As already stated in Section 1, a crystallization $(\Gamma, \gamma)$ of a PL 4-manifold is simple if and only if {\it the 1-skeleton of $K(\Gamma)$ equals the 1-skeleton of a single 4-simplex}. With the notations introduced in the previous section, this is equivalent to require $g_{ijk}=1$ for any distinct $i,j,k \in \Delta_4.$

As a direct consequence of Definition 3, any PL 4-manifold $M$ admitting simple crystallizations turns out to be simply-connected.
On the other hand, any ``standard'' simply-connected PL $4$-manifold (i.e.  $\mathbb S^4$, $\mathbb{CP}^{2}$, $\mathbb{S}^{2} \times \mathbb{S}^{2}$ and the $K3$-surface, together with their connected sums, possibly by taking copies with reversed orientation, too) is proved to admit simple crystallizations (see \cite{[BaS]}).

\bigskip

In the following, we will show that interesting information about simple crystallizations arise by taking into account the handle decompositions induced by the associated coloured triangulations.
First of all, we recall that every closed PL 4-manifold $M$ admits a handle decomposition
$$ M= H^{(0)} \cup (H^{(1)}_1 \cup \dots \cup H^{(1)}_{r_1}) \cup (H^{(2)}_1 \cup \dots \cup H^{(2)}_{r_2}) \cup (H^{(3)}_1 \cup \dots \cup H^{(3)}_{r_3}) \cup  H^{(4)}$$
where $H^{(0)} =\mathbb D^4$ and  each $p$-handle $H^{(p)}_i= \mathbb D^p \times \mathbb D^{4-p}$  ($1 \le p \le 4, \ 1 \le i \le r_p$) is endowed with an an embedding (called {\it attaching map}) $f_i^{(p)}: \partial \mathbb D^p \times \mathbb D^{4-p} \to  \partial (H^{(0)}  \cup \dots (H^{(p-1)}_1 \cup \dots \cup H^{(p-1)}_{r_{p-1}}))$; moreover, it is well-known that $3$- and $4$-handles are attached in a unique way to the union of the $h$-handles, with $0 \leq h \leq 2$.

A standard argument of crystallization theory allows to state that, for any crystallization $(\Gamma,\gamma)$ of a PL 4-manifold $M$ and for any partition
$\{\{i,j,k\},\{r,s\}\}$ of $\Delta_4$,  $M$ admits a decomposition of type $ M \, = \, N(i,j,k) \cup_{\phi} N(r,s),$  where $N(i,j,k)$ (resp. $N(r,s)$) denotes a regular neighbourhood of the subcomplex $K(i,j,k)$ (resp. $K(r,s)$) of $K(\Gamma)$ generated by the vertices labelled $\{i,j,k\}$ (resp. $\{r,s\}$) and $\phi$ is a boundary identification. Any such decomposition turns out to induce a handle decomposition, where $N(i,j,k)$ constitutes the union of the $h$-handles, with $0 \leq h \leq 2$, while $N(r,s)$ is the union of $3$- and $4$-handles.
In the particular case of a simple crystallization, it is easy to prove that a so-called {\it special handlebody decomposition} of $M$ is actually obtained, i.e. a handle decomposition lacking in 1-handles and 3-handles (see \cite[Section 3.3]{[M]}):

\begin{prop} \label{handle-decomposition}
Let $(\Gamma,\gamma)$ be a simple crystallization of a (simply-connected) PL 4-manifold $M$.
Then, for any partition $\{\{i,j,k\},\{r,s\}\}$
of $\Delta_4$, the coloured triangulation $K(\Gamma)$ of $M$ induces a handle decomposition of $M$ consisting of one 0-handle, $g_{rs}-1$ 2-handles and one 4-handle.
\end{prop}

\dimo
Let us fix an arbitrary partition $\{\{i,j,k\},\{r,s\}\}$ of $\Delta_4$. Since $(\Gamma,\gamma)$ is assumed to be a simple crystallization, then $K(r,s)$ consists of exactly one 1-simplex; hence, $N(r,s) \cong_{PL} \mathbb D^4$ trivially follows. On the other hand, the assumption also implies that each of $K(j,k)$, $K(i,k)$ and $K(i,j)$ consists of exactly one 1-simplex; hence, all  $g_{rs}$ 2-simplices of $K(i,j,k)$ have the same boundary.
It is not difficult to check that, if a ``small'' regular neighbourhood of one (arbitrarily fixed) 2-simplex of  $K(i,j,k)$ is considered as a 0-handle $H^{(0)}=\mathbb D^4,$ then the regular neighbourhoods of the remaining $g_{rs}-1$ 2-simplices of  $K(i,j,k)$ may be considered as $g_{rs}-1$ 2-handles attached on its boundary.
Hence,  $N(i,j,k) = H^{(0)} \cup (H^{(2)}_1 \cup \dots \cup H^{(2)}_{g_{rs}-1})$. Moreover, $\partial N(i,j,k) = \partial N(r,s) = \mathbb S^3$.
The proof is completed by noting that the boundary identification $\phi$ between $N(i,j,k)$ and $N(r,s)=\mathbb D^4$ is nothing but the attachment of a 4-handle:
$$ \begin{aligned} M \, = \, N(i,j,k) \cup_{\phi} N(r,s) \,  = & \, [H^{(0)} \cup (H^{(2)}_1 \cup \dots \cup H^{(2)}_{g_{rs}-1})] \cup_{\phi} \mathbb D^4 \, = \\   = \, & H^{(0)} \cup (H^{(2)}_1 \cup \dots \cup H^{(2)}_{g_{rs}-1}) \cup H^{(4)}. \end{aligned} $$
\qed

According to a well-known literature (see, for example, \cite{[M]} or \cite{[K]}), a framed link (possibly with dotted circles) can be associated to any handle decomposition of a PL 4-manifold. If we denote by  $M^3(L,c)$ (resp. $M^4(L,c)$) the 3-manifold (resp. the bounded PL 4-manifold  with boundary $M^3(L,c)$) associated to a framed link $(L, c)$, and - in case $M^3(L,c)$ being PL-homeomorphic to either $\mathbb S^3$ or a connected sum of the orientable $\mathbb S^2$-bundle over $\mathbb S^1$ - we denote by $\bar M^4(L,c)$ the closed  PL 4-manifold associated to  $(L, c),$ then Theorem \ref{framed-link} may be re-stated as follows:

\smallskip
\noindent {\em If $M$ admits simple crystallizations, then a framed link $(L, c)$  (with $\beta_2(M)$ components and no dotted circle) exists, so that
$$ M^3(L,c)= \mathbb S^3  \text{ \ and \ }  \bar M^4(L,c)= M.$$ }

\medskip

\noindent \emph{Proof of Theorem \ref{framed-link}.}

First of all we point out that, if a PL 4-manifold $M$ admits a handle decomposition of the type described in Proposition \ref{handle-decomposition}, then the second Betti number of $M$ must coincide with the number of 2-handles; hence, $g_{rs}-1 = \beta_2(M)$ holds.

Now, let  $(L,c)$ be the framed link obtained by considering, for any $i=1, \dots, \beta_2(M),$  the framed knot in  $\mathbb S^3=\partial H^{(0)}$ corresponding to the attaching map of the $i$-th 2-handle $H^{(2)}_i$ (see \cite[Section 3.1]{[M]}).  The handle decomposition described in Proposition \ref{handle-decomposition} directly ensures  $\bar M^4(L,c)=M$. In particular, no dotted circle appears, since the handle decomposition lacks in 1-handles; moreover, the lacking of 3-handles implies $M^3(L,c)=\mathbb S^3.$
\qed

The following statement  is a straightforward consequence of Proposition \ref{handle-decomposition}: in fact, if an exotic $\mathbb S^4$ or  $\mathbb{CP}^{2}$ exists, its handle decomposition must contain either $1$- or $3$-handles (as pointed out, for example, in \cite{[Ak1]}).

\begin{cor} \label{no-simple}
If an exotic PL-structure on $\mathbb S^4$ (resp. $\mathbb{CP}^{2}$) exists, then the corresponding PL-manifold does not admit simple crystallizations.
\ \  \qed
\end{cor}

\begin{rem}
{\rm In \cite{[C_JKTR]} and \cite{[C_Compl]}, relationships between crystallization theory and (dotted) framed link representation for PL $4$-manifolds are investigated. In particular, a method is described to yield a crystallization $\bar \Lambda(L,c)$ of $\bar M^4 (L,c)$ directly from $(L,c)$. Note that, even if no dotted circle appears and $M^3(L,c)=\mathbb S^3$ is assumed, $\bar \Lambda(L,c)$ is not a (simple) crystallization, and no general procedure is known to obtain a simple crystallization, if any, from $\bar \Lambda(L,c)$. The problem of detecting conditions on $(L,c)$ which ensure the existence of simple crystallizations of $\bar M^4 (L,c)$  could be the matter of a further investigation.}
\end{rem}

\medskip

Let us now state some combinatorial properties of simple crystallizations, which turn out to significantly involve the second Betti number of the represented PL 4-manifold.

\begin{prop} \label{proprieta_cryst_semplici}
Let $(\Gamma,\gamma)$ be an order $2p$ simple crystallization of a (simply-connected) PL 4-manifold $M$.
Then:
\begin{itemize}
\item[(a)] $g_{ij} = 1 + \beta_2(M),$ \ $\forall i,j \in \Delta_4;$
\item[(b)] $p= 1+ 3 \beta_2(M);$
\item[(c)] $ \rho_{\e}(\Gamma)= 2 \beta_2(M),$ \ for any cyclic permutation $\epsilon$ of $\Delta_4.$
\end{itemize}
\end{prop}

\dimo
Within the proof of Theorem \ref{framed-link} we have already noticed that, if a PL 4-manifold $M$ admits a handle decomposition of the type described in Proposition \ref{handle-decomposition},  the second Betti number of $M$ coincides  with the number of 2-handles.
Since such a decomposition, with $g_{rs}-1$ 2-handles, exists for any simple crystallization and for any partition $\{\{i,j,k\},\{r,s\}\}$ of $\Delta_4$,
statement (a)
easily follows.

Let us now apply property (c) of Section 2 to an order $2p$ simple crystallization $(\Gamma,\gamma)$ of a PL 4-manifold $M$ (which - as it is well-known - is simply-connected):
$$\chi (M) = 2 + \beta_2 (M) = - 3p + 10 (1 + \beta_2(M)) - 10 +5,$$
from which $3\beta_2 (M) = p -1$  (i.e. statement (b)) directly follows.

Finally, let us apply equality \eqref{genere} to an order $2p$ simple crystallization $(\Gamma,\gamma)$, by making use of the above statements (a) and (b), too:
$$ 5 (1 + \beta_2(M)) -3 (1 + 3\beta_2(M)) = 2 - 2 \rho_\e(\Gamma),$$
from which
$ \rho_{\e}(\Gamma)= 2 \beta_2(M)$ directly follows, for any cyclic permutation $\epsilon$ of $\Delta_4.$
\qed

Propositions \ref{handle-decomposition} and \ref{proprieta_cryst_semplici} imply that simple crystallizations realize both gem-complexity and regular genus of the represented 4-manifolds and satisfy further combinatorial conditions:

\begin{prop} \label{minimality}
Let $(\Gamma,\gamma)$ be an order $2p$ simple crystallization of a (simply-connected) PL 4-manifold $M$.
Then:
\begin{itemize}
\item[(a)]  $\mathcal G(M)= \rho_{\e}(\Gamma),$ \ for any cyclic permutation $\epsilon$ of $\Delta_4;$
\item[(b)] $k(M)=p-1;$
\item[(c)] $\Gamma$ is rigid and dipole-free.
 \end{itemize}
\end{prop}

\dimo
Statement (a) is a direct consequence of Proposition \ref{proprieta_cryst_semplici}(c), by making use of equation (\ref{genere-beta_2}) of Section 2.

In order to prove statement (b), note that property (d) of Section 2 yields
$$ 2 \sum_{i<j<k} g_{ijk} = 3 \sum_{i<j} g_{ij} -  10p$$ for any order $2p$ crystallization of a PL 4-manifold $M$.  If, further, $M$ is assumed to be simply-connected, by making use of property (c) of Section 2, together with $g_{\hat \imath}=1$ $\forall i \in \Delta_4$ and $\beta_1(M)= \beta_3(M)=0$, the Euler characteristic computation gives
\begin{equation}\label{chi-computation}
2 + \beta_2(M) =   5 - \frac 1 3 \sum_{i<j<k} g_{ijk} + \frac 1 3 p .
\end{equation}
Since $g_{ijk} \ge 1$ trivially holds, we have $ 3 \beta_2(M) \le p-1,$ which proves $k(M) \geq 3 \beta_2(M)$ (already stated as equality (3) of \cite{[CC$_2014$]}).
Let now suppose $(\Gamma,\gamma)$ to be a simple crystallization. In virtue of Proposition \ref{proprieta_cryst_semplici}(a), $k(M) \leq p-1= 3 \beta_2(M)$  follows, too.  Hence, the equality $k(M)=3 \beta_2(M) =p-1$ is established.

Statement (c) is a direct consequence of statement (a), since dipoles and $\rho$-pairs may always be eliminated, yielding another crystallization of $M$, with strictly less order: see \cite[Proposition 7(a)]{[CC$_2014$]}.
\qed

\bigskip
We are now able to prove Theorem \ref{main}, stated in Section 1.

\medskip

\noindent \emph{Proof of Theorem \ref{main}.}

The first statement is a direct consequence of Proposition \ref{proprieta_cryst_semplici}(c) and Proposition \ref{minimality}(a).

\smallskip As far as the second statement is concerned,
note that - in virtue of Proposition \ref{proprieta_cryst_semplici}(b) and Proposition \ref{minimality}(b)
- all PL 4-manifolds admitting a simple crystallization do satisfy condition $k(M) \, = \, 3 \beta_2(M);$ hence, only the reversed implication has to be proved. Then, let $M$ be a simply-connected PL 4-manifolds satisfying $k(M) \, = \, 3 \beta_2(M)$ and let $(\Gamma, \gamma)$ be a crystallization of $M$ realizing its gem-complexity (i.e.: $\#V(\Gamma)= 2(k(M)+1)= 6 \beta_2(M)+2$).
The above equation (\ref{chi-computation}) yields
$$ 2 + \beta_2(M) =   5 - \frac 1 3 \sum_{i<j<k} g_{ijk} + \frac 1 3 (3 \beta_2 (M) +1),$$
and therefore  $$\sum_{i<j<k} g_{ijk} =10$$ directly follows.
Since each summand is at least equal to one,  $g_{ijk} =1$  is proved to hold for any triple $i,j,k \in \Delta_4.$ Hence, $(\Gamma, \gamma)$ turns out to be a simple crystallization, as required.
\qed

\begin{rem}
{\rm Actually, the above proof makes use of the weaker assumption $\beta_1(M)=0$ (instead of the simply-connectedness of $M$), in order to check the existence of a simple crystallization of $M$ when $k(M) = 3 \beta_2(M)$ holds.
Hence, for each orientable PL 4-manifold, the following implication may be stated:

\smallskip
\noindent {\em if $\beta_1(M)=0$ and $ k(M) = 3 \beta_2(M)$, then  $M$ admits a simple crystallization, and hence is simply-connected.} }
\end{rem}

\smallskip

\begin{rem}
{\rm By making use of relations (b) and (c) included in the proof of \cite[Proposition 12]{[CC$_2014$]}, it is not difficult to prove that, if $\pi_1(M)$ is assumed to be trivial, then equality $\mathcal G(M) \, = \, 2 \beta_2(M)$ implies the existence of a crystallization $(\Gamma,\gamma)$ of $M$ and a permutation $\e$ of $\Delta_4$ so that  $\rho_{\e}(\Gamma)= 2 \beta_2(M)$ and $g_{\e_i \e_{i+2} \e_{i+3}}=1$ $\forall i \in \Delta_4.$
However, in general, this does not imply that $(\Gamma,\gamma)$ is simple, since at least one $g_{\e_i \e_{i+1} \e_{i+2}}>1$ may occur. Note that, for example, all rigid dipole-free order $16$ crystallizations satisfy relation $\mathcal G(M) \, = \, 2 \beta_2(M)$, while
 $g_{rst}=2$ for exactly one triple $\{r, s, t\}\subset \Delta_4$ and $g_{ijk}=1$ $\forall \{i,j,k\} \ne\{r, s, t\}.$
}\end{rem}

\bigskip
The characterization of PL 4-manifolds admitting simple crystallizations (Theorem \ref{main}, second statement)
has the following consequence about possible different PL-structures on the same TOP 4-manifold:

\begin{prop} \label{exotic}
Let $M$ and $M^{\prime}$ be two  PL 4-manifolds, with  $M\cong_{TOP} M^{\prime}$ and $M\ncong_{PL} M^{\prime}.$  If both $M$ and $M^{\prime}$ admit a simple crystallization, then $k(M) = k(M^{\prime})$.
\end{prop}

\dimo
It is sufficient to note that $M\cong_{TOP} M^{\prime}$ obviously implies $\beta_2(M)= \beta_2(M^{\prime}),$
and to make use of Proposition \ref{proprieta_cryst_semplici}(b), together with Proposition \ref{minimality}(b).
\qed

\begin{rem}
{\rm In \cite[Section 3]{[CC$_2014$]} an algorithm is described, for the generation of all rigid dipole-free crystallizations of PL 4-manifolds up to a fixed gem-complexity $k$. Proposition \ref{minimality}(c) ensures that such a catalogue must contain all simple crystallizations of PL 4-manifolds whose second Betti number does not exceed $\frac k3$. Actually these catalogues have been generated up to gem-complexity $9$ (\cite{[CC$_2014$]}), hence they present all simple crystallizations of any PL 4-manifold $M$ with $\beta_2(M)\le 3$.
Moreover, Proposition \ref{exotic} guarantees that simple crystallizations representing two distinct PL-structures on the same topological 4-manifold must appear at the same level in the above crystallization catalogues.}
\end{rem}

\bigskip

\section{\hskip -0.7cm . Further results on simple crystallizations}

As already pointed out in Section 3, Basak and Spreer produced a simple crystallization of the $K3$-surface (\cite[Section 7]{[BaS]}), and hence simple crystallizations for any ``standard'' simply-connected PL $4$-manifold are proved to exist.

Theorem \ref{main} has the following consequences about the computation of both PL-invariants
regular genus and gem-complexity for such 4-manifolds:

\begin{prop} \label{calculations}
Let $M \cong_{PL} (\#_r\mathbb {CP}^2)\#(\#_{r^\prime}(-\mathbb {CP}^2)) \# (\#_s(\mathbb S^2\times \mathbb S^2)) \# (\#_t K3) $ \  with $r,r^{\prime},$ $s,t \geq 0$. Then,
$$ \mathcal G(M) = 2(r + r^{\prime} + 2s + 22t)   \quad and \quad k(M) = 3(r + r^{\prime} + 2s + 22t).$$
In particular:
\  $\mathcal G(K3)=44 \ and \ $k(K3)=66$.$
\end{prop}

\dimo
It is sufficient to apply Theorem \ref{main}, by taking into account the values of the second Betti number of each connected summand.
\qed

Moreover, we are able to prove the additivity of both the above invariants under connected sum, within the class of PL 4-manifolds admitting a simple crystallization (and, in particular, for ``standard'' simply-connected PL 4-manifolds).

\begin{prop} \label{additivity}
Let $M$ and $M^{\prime}$ be two (simply-connected) PL 4-manifolds admitting a simple crystallization.
Then:
$$   \mathcal G(M \# M^{\prime}) = \mathcal G(M) + \mathcal G(M^{\prime})  \quad and \quad  k(M \# M^{\prime}) = k(M) + k(M^{\prime}).$$
 \end{prop}

\dimo
It is well-known a general construction - called {\it graph-connected sum} - yielding, from any gem $\Gamma$  (resp. $\Gamma^{\prime}$) of the PL $n$-manifold $M$  (resp. $M^{\prime}$), a gem $\Gamma \# \Gamma^{\prime}$ of $M \# M^{\prime}$: see \cite{[FGG]} for details.

On the other hand, it is not difficult to check that, if both $\Gamma$ and $\Gamma^{\prime}$ are simple crystallizations, then $\Gamma \# \Gamma^{\prime}$ is, too.
Hence, since simple crystallizations do always realize both regular genus and gem-complexity
of the represented PL 4-manifolds (Proposition  \ref{minimality}(a) and (b)), the thesis easily follows.
\vspace{-0.3truecm} \ \ \qed

\begin{rem}
{\rm Note that the relation $\mathcal G(M \# M^{\prime}) \le  \mathcal G(M) + \mathcal G(M^{\prime})$ can be stated for all PL $n$-manifolds by direct estimation of $\mathcal G(M \# M^{\prime})$
on any gem $\Gamma \# \Gamma^{\prime},$ when $\Gamma,$ $\Gamma^{\prime}$ are assumed to be gems of $M, M^{\prime}$ realizing regular genus of the represented $n$-manifolds. Moreover - as pointed out in Section 1 - the additivity of regular genus under connected sum has been conjectured, and the associated (open) problem is significant especially in dimension four.

In \cite[Corollary 4]{[GG]}, two classes of closed (not necessarily orientable) $4$-manifolds have been detected, for which additivity of regular genus holds. It is not difficult to check that the first one (characterized by relation $\mathcal G(M)= 2 \chi(M) -4$) includes - in virtue of Theorem \ref{main} - all 4-manifolds admitting a simple crystallization, while the second one (characterized by relation $\mathcal G(M)=1 - \frac {\chi(M)} 2 $) consists - in virtue of \cite[Proposition 2]{Casali-Malagoli} - of connected sums of $\mathbb S^3$-bundles over $\mathbb S^1.$
}
\end{rem}

\smallskip

We conclude the paper by reporting two results already proved in \cite{[CC$_2014$]}. The first one concerns the existence of simple crystallizations of standard simply-connected PL 4-manifolds with $\beta_2 \leq 2 $, and has been obtained as a direct consequence of 4-dimensional crystallization catalogues:

\begin{prop}  {\rm (\cite[Proposition 17]{[CC$_2014$]})}
\ \par
\begin{itemize}
\item{}  $\mathbb S^4$ and $\mathbb{CP}^{2}$ admit a unique simple crystallization;
\item{}  $\mathbb{S}^{2} \times \mathbb{S}^{2}$ admits exactly $267$ simple crystallizations;
\item{}  $\mathbb{CP}^{2} \# \mathbb{CP}^{2}$ admits exactly $583$ simple crystallizations;
\item{}  $\mathbb{CP}^{2} \# (- \mathbb{CP}^{2})$ admits exactly $258$ simple crystallizations.
    \end{itemize}
\vspace{-0.3truecm} \rightline{$\Box$ \ \ \ \ \ \ \ }
\end{prop}

Finally, the existence of simple crystallizations can be related to known results and open problems about exotic structures on ``standard'' simply-connected PL $4$-manifolds
(see Proposition \ref{exotic}):

\begin{prop} {\rm (\cite[Proposition 18]{[CC$_2014$]})} \label{exotic-simple}
\ \par
\begin{itemize}
\item[(a)] Let $M$ be  $\mathbb S^4$ or $\mathbb{CP}^{2}$ or $\mathbb{S}^{2} \times \mathbb{S}^{2}$ or $\mathbb{CP}^{2} \# \mathbb{CP}^{2}$ or $\mathbb{CP}^{2} \# (- \mathbb{CP}^{2})$; if an exotic PL-structure on $M$ exists, then the corresponding PL-manifold does not admit a simple crystallization.
\item[(b)] Let $\bar M$ be a PL $4$-manifold TOP-homeomorphic but not PL-homeomorphic to $\mathbb{CP}^{2} \#_2 (-\mathbb{CP}^{2})$; then, either $\bar M$ does not admit a simple crystallization, or $\bar M$ admits an order $20$ simple crystallization (i.e.: $k(\bar M)= 9 = k (\mathbb{CP}^{2} \#_2 (-\mathbb{CP}^{2}))$).
\item[(c)]  Let $r \in \{3,5,7,9,11,13 \} \cup \{r =4n-1 \ / \ n \ge 4\} \cup \{r =4n-2 \ / \ n \ge 23\}$; then, infinitely many simply-connected PL $4$-manifolds with $\beta_2=r$ do not admit a simple crystallization.
\end{itemize}
\vspace{-0.3truecm} \rightline{$\Box$ \ \ \ \ \ \ \ }
\end{prop}

\begin{rem}
{\rm Note that, while in \cite{[CC$_2014$]} the proof of Proposition \ref{exotic-simple}(a) is directly based on the analysis of the crystallization catalogue up to gem-complexity eight, in the present paper the cases regarding $\mathbb S^4$ and $\mathbb{CP}^{2}$ have already been obtained via Proposition \ref{handle-decomposition}: see Corollary \ref{no-simple}.}
\end{rem}

\bigskip

\par \noindent {\bf Acknowledgements.} This work was supported by the ``National Group for Algebraic and Geometric Structures, and their Applications'' (GNSAGA - INDAM) and by M.I.U.R. of Italy (project ``Strutture Geometriche, Combinatoria e loro Applicazioni'').

\bigskip

{\footnotesize{
}}

\begin{thebibliography}{survey}

\bibitem{[Ak1]} S. Akbulut, {\em The Dolgachev Surface  - Disproving Harer-Kas-Kirby Conjecture}
Comm. Math. Helvetici {\bf 87(1)} (2012), 187-241.


\bibitem{[Ak2]} S. Akbulut, {\em An infinite family of exotic Dolgachev surfaces without 1- and 3- handles},
Jour. of GGT {\bf 3} (2009), 22-43.


\bibitem{survey} P. Bandieri - M. R. Casali - P. Cristofori - L. Grasselli - M. Mulazzani, {\em Computational aspects of crystallization theory: complexity, catalogues  and classification of 3-manifolds}, Atti Sem. Mat. Fis. Univ. Modena \textbf{58} (2011), 11--45.

\bibitem{[BCrG$_1$]}
P. Bandieri - P. Cristofori - C. Gagliardi, {\em Nonorientable 3-manifolds admitting coloured triangulations with at most 30 tetrahedra}, J. Knot Theory Ramifications {\bf 18} (2009), 381-395.

\bibitem{[BaD]} B. Basak - B. Datta, {\em Minimal crystallizations of 3-manifolds}, Electron. J. Combinatorics {\bf 21}(1) (2014), \#P1.61, 1-25.

\bibitem{[BaS]} B. Basak - J. Spreer, {\em Simple crystallizations of 4-manifolds}, Advances in Geometry (to appear), arXiv:1407.0752.

\bibitem{[C$_1992$]} M. R. Casali, {\em A combinatorial characterization of 4-dimensional handlebodies},  Forum Math.  {\bf  4} (1992), 123-134.

\bibitem{[Casali_UMI]}	M. R. Casali, {\em An infinite class of bounded 4-manifolds having regular genus three}, Bollettino Un. Mat. Ital. {\bf 10-A} (1996), 279-303.

\bibitem{Casali} M. R. Casali, {\em Classifying PL 5-manifolds by regular genus: the boundary  case}, Canadian J. Math. {\bf 49} (1997), 193--211.

\bibitem{[C$_2$]} M. R. Casali, {\em Classification of non-orientable 3-manifolds admitting decompositions into $\le 26$ coloured tetrahedra}, Acta Appl. Math. {\bf 54} (1999), 75-97.

\bibitem{[C_JKTR]} 	M.R.Casali,  {\em From framed links to crystallizations of bounded 4-manifolds}, Journal of Knot Theory and its Ramifications {\bf 9(4)} (2000), 443-458.

\bibitem{[C_Compl]} M.R.Casali, {\em Dotted links, Heegaard diagrams and coloured graphs for PL 4-manifolds}, Revista Matematica Complutense {\bf 17(2)} (2004), 435-457.

\bibitem{[C$_4$]} M. R. Casali, {\em Computing Matveev's complexity of non-orientable 3-manifolds via crystallization theory}, Topology Appl. {\bf 144} (2004), 201-209.

\bibitem{Casali-Oberwolfach} M. R. Casali, {\em Catalogues of PL-manifolds and complexity estimations via crystallization theory}, Oberwolfach Report no. \textbf{24}/{2012}, 58--61.  [DOI: 10.4171/OWR/2012/24]

\bibitem{[CC$_1$]} M. R. Casali - P. Cristofori, {\em Computing Matveev's complexity via crystallization theory: the orientable case}, Acta Appl. Math. 92 (2) (2006), 113-123.

\bibitem{[CC$_2$]} M. R. Casali - P. Cristofori, {\em A catalogue of orientable 3-manifolds triangulated by $30$ coloured tetrahedra}, J. Knot Theory Ramifications {\bf 17} (2008), 1-23.

\bibitem{[CC$_3$]} M. R. Casali - P. Cristofori: {\em A note about complexity of lens spaces}. Forum Math. (2014), DOI: 10.1515/forum-2013-0185, published online February 19, 2014. [arXiv:1309.5728]

\bibitem{[CC$_2014$]}  M. R. Casali - P. Cristofori, {\em Cataloguing PL 4-manifolds by gem-complexity}, arXiv:1408.0378.

\bibitem{[CG]}  M. R. Casali - C. Gagliardi, {\em Classifying PL 5-manifolds up to regular genus seven}, Proc. Amer. Math. Soc. {\bf 120} (1994), 275-283.

\bibitem{Casali-Malagoli} M. R. Casali - L. Malagoli,  {\em Handle-decompositions of PL 4-manifolds}, Cahiers Topologie Geom. Differentielle Categ. {\bf 38} (1997), 141--160.

\bibitem{[Cav$_1992$]} A. Cavicchioli, {\em On the genus of smooth 4-manifolds},  Trans. Amer. Math. Soc.  {\bf  331 (1)} (1992), 203-214.

\bibitem{[FG$_2$]}  M. Ferri -  C. Gagliardi, {\em The only genus zero -manifold is $\mathbb S^n$}, Proc. Amer. Math. Soc. {\bf 85} (1982), 638-642.

\bibitem{[FGG]} M. Ferri - C. Gagliardi - L.Grasselli, {\em A graph-theoretical representation of PL-manifolds. A survey on crystallizations},   Aequationes Math. {\bf 31} (1986), 121-141.

\bibitem{[G]} C. Gagliardi, {\em Extending the concept of genus to dimension $n$}, Proc. Amer. Math. Soc. {\bf 81} (1981), 473-481.

\bibitem{[GG]} C. Gagliardi - L.Grasselli, {\em Representing products of polyhedra by products of edge-coloured graphs}, Journal of Graph Theory {\bf 17} (1993), 549-579.

\bibitem{[K]} R. Kirby, {\em The topology of 4-manifolds}, Lecture Notes in Math. {\bf 1374}, Springer-Verlag, 1989.

\bibitem{[L]} S. Lins, {\em Gems, computers and attractors for 3-manifolds}, Knots and Everything {\bf 5}, World Scientific, 1995.

\bibitem{[M]} R. Mandelbaum, {\em Four-dimensional topology: an introduction}, Bull. Amer. Math. Soc. {\bf 2} (1980), 1-159.

\bibitem{[Mat]} S.Matveev, {\em Complexity theory of three-dimensional manifolds}, Acta Appl. Math. {\bf 19} (1990),
101-130.

\bibitem{[SK]} J. Spreer - W. Kuhnel, {\em  Combinatorial properties of the K3 surface: Simplicial blowups and slicings}, Experiment. Math. {\bf 20(2)} (2011), 201-216.

\bibitem{[Sw]} E. Swartz, {\em The average dual surface of a cohomology class and minimal
simplicial decompositions of infinitely many lens spaces}, arXiv:1310.1991.




\end{thebibliography}
\end{document}